\documentclass[11pt]{article} 
\usepackage{amsmath}
\usepackage{amsfonts}
\usepackage{amssymb}
\usepackage{graphicx}
\usepackage{color}
\usepackage[latin1]{inputenc}
\usepackage[T1]{fontenc}
\usepackage[french]{babel}

\makeatletter
\def\@sect#1#2#3#4#5#6[#7]#8{%
  \ifnum #2>\c@secnumdepth
    \let\@svsec\@empty
  \else
    \refstepcounter{#1}%
    \protected@edef\@svsec{\@seccntformat{#1}\relax}%
  \fi
  \@tempskipa #5\relax
  \ifdim \@tempskipa>\z@
    \begingroup
      #6{%
        \@hangfrom{\hskip #3\relax\@svsec}%
          \interlinepenalty \@M #8\@@par}%
    \endgroup
    \csname #1mark\endcsname{#7}%
    \addcontentsline{toc}{#1}{%
      \ifnum #2>\c@secnumdepth \else
        \protect\numberline{\csname the#1\endcsname.}%
      \fi
      #7}%
  \else
    \def\@svsechd{%
      #6{\hskip #3\relax
      \@svsec #8}%
      \csname #1mark\endcsname{#7}%
      \addcontentsline{toc}{#1}{%
        \ifnum #2>\c@secnumdepth \else
          \protect\numberline{\csname the#1\endcsname.}%
        \fi
        #7}}%
  \fi
  \@xsect{#5}}
\def\@seccntformat#1{\csname the#1\endcsname.\quad}
\makeatother
\newtheorem{theo}[equation]{Th\'eor\`eme}

\newtheorem{proposition}[equation]{Proposition}
\newtheorem{df}[equation]{D\'efinition}
\newtheorem{cor}[equation]{Corollaire}
\newtheorem{question}[equation]{Question}

\newenvironment{remarque}{
\refstepcounter{equation}\trivlist%
\item[\hskip \labelsep{\bfseries Remarque \theequation.\ }]}%
{\endtrivlist}%

\makeatletter
\renewcommand\theequation{\thesection.\arabic{equation}}
\@addtoreset{equation}{section}
\makeatother
\newcommand{\carrenoir}{\rule{0.5em}{0.5em}}


\newcommand{\oper}[2]{\newcommand{#1}{\mathop{\mathrm{#2}}\nolimits} }
\oper{\Vol}{Vol}
\newcommand{\R}{\mathbb R}
\newcommand{\de}{\mathrm{ d }}
\oper{\scal}{scal}

\DeclareSymbolFont{greek}{OML}{ptmcm}{m}{it}
\DeclareMathSymbol{\codiff}{\mathord}{greek}{"0E}

\newcommand{\N}{\mathbb N}

\title{Extrema de valeurs propres dans une classe conforme}
\author{Pierre Jammes}
\date{}
\begin{document}
\maketitle
{\small 
\textsc{Résumé.---}
On s'intéresse au problème de savoir quelle est la rigidité apportée au 
spectre d'une variété riemannienne compacte par le fait de fixer son volume 
et se classe conforme, et en particulier de déterminer si on peut faire
tendre les valeurs propres vers 0 ou l'infini sous cette contrainte.
On considère successivement les cas du laplacien usuel agissant sur
les fonctions, l'opérateur de Dirac, le laplacien conforme et le
laplacien de Hodge-de~Rham.

Mots-clefs : valeurs propres, géométrie conforme, métriques extrémales.

MSC2000 : 35P15, 58J50, 58E11}

\section{Introduction}
Soit $M$ une variété compacte de dimension $n\geq2$. Si $M$ est
munie d'une métrique riemannienne $g$, on peut définir le 
laplacien $\Delta=-\mathrm{div}\ \mathrm{grad}$, agissant sur les
fonctions de $M$, dont nous noterons le spectre 
$0=\lambda_0(M,g)<\lambda_1(M,g)\leq\lambda_2(M,g)\leq\ldots$

Un problème classique consiste à étudier la fonctionnelle 
$g\mapsto\lambda_k(M,g)$ sous certaines contraintes géométriques
et d'en chercher les métriques extrémales. 
Ce problème peut se généraliser à d'autres opérateurs.

 Nous nous intéresserons ici à l'étude de la fonctionnelle 
\begin{equation}\label{intro:eq}
g\mapsto\lambda_k(M,g)\Vol(M,g)^{\frac2n},
\end{equation}
le facteur volume permettant de normaliser cette fonctionnelle en la 
rendant invariante par homothétie, en nous restreignant à l'espace des 
métriques riemanniennes conformes à une métrique $g$ donnée quelconque, 
c'est-à-dire à la classe conforme $[g]=\{h^2g,\ h\in C^\infty(M),\ h>0\}$, 
ce sujet ayant connu récemment un grand nombre de développements. 
 
 On verra successivement le cas du laplacien usuel agissant sur les
fonctions (section~\ref{fonctions}), l'opérateur de Dirac ---~ou plus 
précisément son carré pour que la fonctionnelles (\ref{intro:eq}) reste
invariante par homothétie~--- et le laplacien conforme dans la 
section~\ref{dirac}, et le laplacien de Hodge-de~Rham dans la 
section~\ref{formes}.

 Je remercie Bernd Ammann et Bruno Colbois pour de nombreuses discussions
autour de ce sujet, ainsi que Gérard Besson qui 
m'a permis d'exposer mes travaux au séminaire de théorie spectrale et 
géométrie de l'institut Fourier.
\section{Laplacien usuel}\label{fonctions}
\subsection{Majoration des valeurs propres et spectre conforme}
Si $(M,g)$ est une variété riemannienne compacte, nous noterons le spectre 
du laplacien usuel $\Delta:f\mapsto \nabla^*\nabla f$, agissant sur 
les fonctions de $M$, par 
\begin{equation}
0=\lambda_0(M,g)<\lambda_1(M,g)\leq\lambda_2(M,g)\leq\ldots
\end{equation}
À volume fixé, on peut facilement faire tendre ces valeurs propres vers
zéro en construisant des haltères de Cheeger, et on peut aussi les faire 
tendre vers l'infini d'après \cite{cd94}. Si on fixe la classe conforme,
on peut encore exhiber des petites valeurs propres, mais en revanche
N.~Korevaar a montré que pour tout $k$, $\lambda_k(M,g)$ est uniformément
majorée:
\begin{theo}[\cite{ko93}]
Soit $M$ une variété compacte de dimension $n$ et $C$ une classe conforme
de métriques sur $M$. Il existe une constante $a(C)>0$ telle que
$$\sup_{g\in C}\lambda_k(M,g)\Vol(M,g)^{\frac2n}\leq a\cdot k^{\frac2n}.$$
\end{theo}
Pour $k=1$, une telle majoration avait déjà été obtenue par A.~El~Soufi et
S.~Ilias dans \cite{esi86} (voir paragraphe suivant).

Ce résultat permet de définir des invariants conformes de la variété
$M$ en considérant la borne supérieure de $\lambda_k(M,g)\Vol(M,g)^{\frac2n}$
dans une classe conforme. ces invariants ont été définis et étudiés par
B.~Colbois et A.~El~Soufi dans \cite{ces03} :
\begin{df}[\cite{ces03}]
Soit $M$ une variété riemannienne compacte et $C$ une classe conforme
de métriques sur $M$. La $k$-ième valeur propre conforme de $(M,C)$ est
définie par
$$\lambda_k^c(M,C)=\sup_{g\in C}\lambda_k(M,g)\Vol(M,g)^{\frac2n}.$$
La suite $(\lambda_k^c(M,C))_k$ est appelée spectre conforme de $(M,C)$.
\end{df}
On connaît peu de valeurs exactes des $\lambda_k^c(M,C)$. Outre le
cas de la 1\iere{} valeur propre de certaines variétés
s'immergeant isométriquement minimalement dans des sphères que nous
verrons au paragraphe~\ref{fonctions:vc}, le seul exemple semble être
$\lambda_2^c(S^2,C_\mathrm{can})=16\pi$ calculé par N.~Nadirashvili 
dans \cite{na02}.

B.~Colbois et A.~El~Soufi obtiennent dans \cite{ces03} quelques
estimées générales des valeurs propres conformes, à commencer par
une minoration par le spectre conforme de la sphère:
\begin{theo}
Pour toute classe conforme $C$ sur une variété $M$ de dimension $n$, et
tout entier $k>0$ on a
$$\lambda_k^c(M,C)\geq\lambda_k^c(S^n,C_\mathrm{can}).$$
\end{theo}
Ce résultat était déjà connu pour $k=1$ (\cite{fn99}).

Leur second résultat est que la différence entre deux valeurs propres
conformes consécutives est uniformément minorée :
\begin{theo}\label{fonctions:th3}
Pour toute classe conforme $C$ sur une variété $M$ de dimension $n$, et
tout entier $k>0$ on a
$$\lambda_{k+1}^c(M,C)^{\frac n2}-\lambda_k^c(M,C)^{\frac n2}\geq
\lambda_1^c(S^n,C_\mathrm{can})^{\frac n2}=n^{\frac n2}\omega_n,$$
où $\omega_n$ désigne le volume de $S^n$ pour sa métrique canonique.
\end{theo}
On peut en déduire une minoration de $\lambda_k^c(M,C)$ en fonction de
$n$ et $k$ :
\begin{cor}
Pour toute classe conforme $C$ sur une variété $M$ de dimension $n$, et
tout entier $k>0$, on a
$$\lambda_k^c(M,C)\geq n\omega_n^{\frac2n}k^{\frac2n}.$$
\end{cor}
Notons que dans cette dernière égalité, on a égalité sur les sphères
pour $k=1$ et $n$ quelconque, ainsi que pour $k=n=2$.

La question se pose naturellement de savoir quelles métriques
peuvent réaliser le maximum de $\lambda_k(M,g)$ dans une classe
conforme. A.~El~Soufi et S.~Ilias montre dans \cite{esi03} que
pour une telle métrique, la valeur propre $\lambda_k(M,g)$ est
nécessairement multiple. On peut en déduire, en conjonction avec
le théorème~\ref{fonctions:th3} qu'une métrique ne peut pas
maximiser trois valeurs propres consécutives dans sa classe conforme.

\subsection{Volume conforme et 1\iere{} valeur propre conforme}%
\label{fonctions:vc}
 Le cas de la première valeur propre du laplacien est particulier car
A.~El~Soufi et S.~Ilias ont donné une majoration plus précise que
celle de N.~Korevaar ---~et en un sens optimale~--- en utilisant le volume 
conforme. 

Si on note $G_N$ le groupe des difféomorphismes conforme de la sphère
$S^N$ munie de sa structure conforme canonique, le volume conforme
de la variété $M$ pour la classe conforme $C$, introduit par
P.~Li et S.~T.~Yau dans \cite{ly82}, est défini par
\begin{equation}\label{fonctions:defvc}
V_\mathrm{c}(M,C)=\inf_N\inf_{\varphi:(M,C)\to S^N}\sup_{\gamma\in G_N}
\Vol(\gamma\circ\varphi(M)),
\end{equation}
l'application $\varphi$ parcourant l'ensemble des difféomorphismes
conformes de $(M,C)$ dans $S^N$.

La propriété du volume conforme qui nous intéresse ici est la suivante :
\begin{theo}\label{fonctions:vcth1}
Si $(M,g)$ est une variété riemannienne compacte, alors
$$\lambda_1(M,g)\Vol(M,g)^{\frac2n}\leq nV_\mathrm{c}(M,[g])^{\frac2n}.$$
\end{theo}
Ce théorème est démontré pour $n=2$ dans \cite{ly82}, et généralisé
en toute dimension par A.~El~Soufi et S.~Ilias dans \cite{esi86}.

La principale motivation de P.~Li et S.~T.~Yau était d'utiliser le 
théorème~\ref{fonctions:vcth1} pour démontrer des cas particuliers de 
la conjecture de Willmore. Signalons que cette inégalité a trouvé
récemment une nouvelle application : en la généralisant à des
orbivariétés, I.~Agol a montré dans \cite{ag05} qu'il n'y a qu'un
nombre fini de groupes de réflexions kleiniens arithmétiques maximaux.

L'inégalité du théorème~\ref{fonctions:vcth1} a l'intérêt qu'on peut
caractériser le cas d'égalité et en exhiber des exemples. 
\begin{df}[\cite{esi92}]
On dit qu'une métrique $g$ sur une variété compacte $M$ est 
$\lambda_1$-minimale si $(M,g)$ s'immerge isométriquement minimalement
dans une sphère par ses premières fonctions propres.
\end{df}
\begin{theo}[\cite{esi86}]
Soit $(M,g)$ est une variété riemannienne compacte. On a l'égalité
$\lambda_1(M,g)\Vol(M,g)^{\frac2n}=nV_\mathrm{c}(M,[g])^{\frac2n}$
si et seulement si $(M,g)$ est $\lambda_1$-minimale.
\end{theo}
\begin{remarque}
La $\lambda_1$-minimalité est aussi fortement liée à l'extrémalité de la 
métrique pour la fonctionnelle (\ref{intro:eq}) dans l'espace de toute 
les métriques (voir \cite{na96} et \cite{esi00}).
\end{remarque}
 La condition de $\lambda_1$-minimalité a été étudiée par T.~Takahashi dans
\cite{ta66}, où il montre qu'elle est équivalente à l'existence d'une
famille de premières fonctions propres $(f_1,\ldots,f_k)$ telle que
$g=\sum_{i=1}^k\de f_i^2$, et que cette condition est vérifiée sur les
variétés homogènes irréductibles. Elle l'est aussi sur les variétés
fortement harmoniques (\cite{be78})
et sur quelques autres exemples (\cite{mou84}). On peut en déduire les
valeurs exactes du volume conforme et de la première valeur propre
conforme sur les espaces projectifs munis de leur métrique
canonique :
\renewcommand{\arraystretch}{1.5}
\begin{center}
\begin{tabular}{|c|c|c|c|}
\hline
variété&dimension&$V_{\mathrm c}$&$\lambda_1^c$\\\hline
$S^n$ & $n$ & $\omega_n$ & $n\omega_n^{\frac2n}$\\
$\R P^n$ & $n$ & $\left(\frac{2(n+1)}n\right)^{\frac n2}\frac{\omega_n}2$ &
$2^{\frac{n-2}n}(n+1)\omega_n^{\frac2n}$\\
$\mathbb CP^n$ & $2n$ & $\left(\frac{2\pi(n+1)}n\right)^n\frac1{n!}$ &
$4\pi(n+1){n!}^{-\frac1n}$\\
$\mathbb HP^n$ & $4n$ & $\left(\frac{2\pi(n+1)}n\right)^{2n}\frac1{(2n+1)!}$ &
$8\pi(n+1)(2n+1)!^{-\frac1{2n}}$\\
$\mathbb C\mathrm aP^2$ & $16$ & $\frac{6(3\pi)^8}{11!}$ & 
$48\pi\left(\frac6{11!}\right)^{\frac18}$
\\\hline
\end{tabular}
\end{center}
\renewcommand{\arraystretch}{1}
\begin{remarque}
En dimension~2, la fonctionnelle (\ref{intro:eq}) est uniformément
majorée sur l'espace de toute les métriques, et une métrique réalisant 
son maximum est $\lambda_1$-minimale. Outre le cas de la sphère, on 
connaît deux valeurs explicites de $\lambda_1^c$ dans cette situation.
Sur le tore, le maximum est atteint par le tore équilatéral 
$T^2_{\mathrm{eq}}=\R^2/\left((1,0)\mathbb Z\oplus(\frac12,\frac{\sqrt3}2)
\mathbb Z\right)$ et $\lambda_1^c(T_{\mathrm{eq}})=\frac{8\sqrt3}3\pi^2$
(voir \cite{na96}). 
Ce maximum a aussi été récemment calculé sur la bouteille de Klein
dans \cite{esgj07}.
\end{remarque}
\begin{remarque}
La propriété de $\lambda_1$-minimalité est stable par produit:
on peut facilement vérifier que si deux variétés ont même première valeur 
propre et que leur métrique peut s'écrire sous la forme 
$g=\sum_{i=1}^k\de f_i^2$ où les $f_i$ sont des premières fonctions propres,
il en va de même pour leur produit. On peut construire ainsi d'autres
exemples comme le tore de Clifford $S^1\times S^1$.
\end{remarque}

\subsection{Spectre conforme minimal}\label{fonctions:scm}
L.~Friedlander et N.~Nadirashvili ont redémontré de manière indépendante 
dans \cite{fn99} que $\lambda_1(M,g)\Vol(M,g)^{\frac2n}$ est majoré
dans une classe conforme. Ils définissent aussi un nouvel invariant
différentiable de $M$ par $\nu(M)=\inf_g\lambda_1^c(M,[g])$ et montrent qu'il
est toujours minoré par $\lambda_1^c(S^n,C_{\text{can}})$. On peut
définir un «~spectre conforme minimal~» pour $M$ en étendant cet invariant
aux autres valeurs propres par 
\begin{equation}
\nu_k(M)=\inf_g\lambda_k^c(M,[g]).
\end{equation}
Ces invariants sont très mal connus. On peut déduire des résultats évoqués
précédemment que $\nu_1(S^n)=\lambda_1^c(S^n,C_{\text{can}})$,
on peut aussi affirmer que $\nu_2(S^2)=16\pi$ et 
$\nu_1(\R P^2)=12\pi$ en utilisant le fait que ces deux variétés 
n'ont qu'une classe conforme. A.~Girouard a récemment montré dans 
\cite{gi05} que $\nu_1(T^2)=\nu_1(K^2)=8\pi$ où $K^2$ désigne la bouteille
de Klein, et on conjecture que pour toute surface compacte $\Sigma$ autre 
que $\R P^2$, on a $\nu_1(\Sigma)=8\pi$.

Il est difficile d'extrapoler ce que peut être le
comportement de ces invariants à partir d'aussi peu d'exemples, surtout
que la dimension~2 pourrait être pathologique. Le cas du plan
projectif montre cependant que l'invariant $\nu_1(M)$ n'est pas
trivial.

Un approche possible pour obtenir des estimées de $\nu_1(M)$ est 
d'utiliser le volume conforme pour définir un nouvel invariant différentiel:
\begin{df}
Si $M$ est une variété différentielle compacte, on définit le \emph{volume
conforme minimal de $M$} par
$$V_\mathrm{cm}(M)=\inf_gV_\mathrm{c}(M,[g]).$$
\end{df}
On a alors la majoration suivante, qui découle immédiatement du 
théorème~\ref{fonctions:vcth1}:
\begin{equation}
\nu_1(M)\leq nV_\mathrm{cm}(M)^{\frac2n}.
\end{equation}
\begin{remarque}
Pour définir le volume conforme minimal, il suffit dans la
relation (\ref{fonctions:defvc}) de supprimer la condition que 
le difféomorphisme $\varphi$ est conforme. Contrairement au spectre
conforme minimal dont la définition fait appel à des
structures riemanniennes et conformes sur la variété, on a seulement 
besoin ici de sa structure différentielle.
\end{remarque}

Les seules valeurs exactes connues du volume conforme minimal sont
$V_\mathrm{cm}(S^n)=\omega_n$ et $V_\mathrm{cm}(\R P^2)=6\pi$. 
On peut cependant montrer qu'à dimension fixée, il est uniformément
majoré:
\begin{theo}[\cite{ja08}]
Pour toute entier $n\geq2$, il existe une constante $c(n)>0$ telle que
pour toute variété compacte $M$ de dimension $n$, on a $V_\mathrm{cm}(M)\leq
c(n)$.
\end{theo}
Le principe de la démonstration est d'étudier comment varie le
volume conforme minimal quand on pratique des chirurgies sur la variété.
\section{Opérateur de Dirac et laplacien conforme}\label{dirac}
\subsection{Notations}
 Je rassemble ici deux opérateurs qui semblent assez différent par
construction mais qui se révèlent avoir de nombreuses propriétés
communes.

Si $(M,g)$ est une variété riemannienne munie d'une structure spinorielle
$\chi$, on peut définir l'opérateur de Dirac $\mathrm D$ agissant sur les 
sections du fibré des spineurs $\Sigma_gM$ (voir par exemple \cite{hi01}). 
C'est un opérateur 
elliptique du 1\ier{} ordre dont le spectre n'est pas borné inférieurement. 
Par commodité, et pour garder une certaine cohérence des notations par rapport 
aux autres opérateurs traités dans ce texte, nous considérerons le
spectre de $\mathrm D^2$, que nous noterons
\begin{equation}
0\leq\lambda_1(M,\chi,g)\leq\lambda_2(M,\chi,g)\leq\ldots.
\end{equation}
On notera en outre $\lambda^+(M,\chi,g)$ la première valeur propre
strictement positive.

Le laplacien conforme, appelé aussi opérateur de Yamabe, qui agit sur les
fonctions de $M$ est défini par
\begin{equation}
\mathrm L_g:f\mapsto\Delta f+\frac{n-2}{4(n-1)}\scal_g f,
\end{equation}
où $\scal_g$ désigne la courbure scalaire pour la métrique $g$.
Le spectre de $\mathrm L$ est borné inférieurement, mais peut contenir un 
nombre fini de valeur propres négatives. Nous noterons son spectre 
complet par
\begin{equation}
\lambda_1(M,\mathrm L_g)\leq\lambda_2(M,\mathrm L_g)\leq\ldots
\end{equation}
et par $\lambda^+(M,\mathrm L_g)$ (resp. $\lambda^-(M,\mathrm L_g)$) la plus 
petite valeur propre (en valeur absolue) strictement positive (resp. 
strictement négative).

Les opérateurs $\mathrm D$ et $\mathrm L$ ont en commun une propriété de 
covariance conforme: on dit qu'on opérateur $\mathrm T$ d'ordre $j$ est 
conformément covariant s'il vérifie la relation
\begin{equation}
\mathrm T_{h^2g}=h^{-\frac{n+j}2}\mathrm T_gh^{\frac{n-j}2}
\end{equation}
où $h$ est une fonction strictement positive
et $\mathrm T_g$ désigne l'opérateur $\mathrm T$ pour la métrique $g$. 
En particulier, la dimension des noyaux de $\mathrm D$ et $\mathrm L$ sont 
des invariants conformes, ce qui justifie le fait de ne s'intéresser ici
qu'à leurs valeurs propres non nulles.

Certaines connivences entre les spectres de ces deux opérateurs sont bien 
connues, citons par exemple l'inégalité suivante due à O.~Hijazi 
(\cite{hi86}) sur laquelle nous reviendront:
\begin{equation}\label{dirac:hijazi}
\lambda_1(M,\chi,g)\geq\frac n{4(n-1)}\lambda_1(M,\mathrm L_g).
\end{equation}
Cette inégalité n'est significative que si $\lambda_1(M,\mathrm L_g)$ est 
strictement positif.
\subsection{Minoration conforme du spectre}
Le comportement du spectre de ces deux opérateurs contraste avec le cas
du laplacien usuel. Le premier résultat obtenu est qu'on ne peut pas
faire tendre les valeurs propres vers zéro dans une classe conforme:
\begin{theo}\label{dirac:min}
Soit $(M,\chi)$ une variété spinorielle compacte et $C$ une classe conforme 
de métrique sur $M$. Alors
$$\inf_{g\in C}\lambda^+(M,\chi,g)\Vol(M,g)^\frac2n>0.$$
\end{theo}
Ce théorème à été démontré par J.~Lott (\cite{lo86}) dans le cas
ou l'opérateur de Dirac est inversible, et dans le cas général par B.~Ammann
(\cite{am03}). Les deux auteurs remarquent que la démonstration
est valable pour n'importe quel opérateur elliptique autoadjoint 
conformément covariant de degré strictement inférieur à $n$, et quel que 
soit le signe des valeurs propres. 
On peut donc remplacer $\lambda^+(M,\chi,g)$ par $\lambda^+(M,\mathrm L_g)$
ou $|\lambda^-(M,\mathrm L_g)|$ dans le théorème~\ref{dirac:min}.

Dans \cite{am03b}, B.~Ammann étudie l'existence de métriques réalisant
la borne inférieure du théorème~\ref{dirac:min}. Dans le cas où cette 
borne est plus petite que celle de la sphère canonique (cette condition 
est étudiée plus en détail dans \cite{ahm06}, elle est en particulier 
vérifiée par $M=\R P^{4k+3}$), il montre 
qu'elle est atteinte en élargissant la classe conforme $C$ à des métriques
dégénérées de la forme $h^2g$ où $g$ est une métrique de $C$ et
$h$ une fonction $C^{2,\alpha}$ qui peut s'annuler. 

 Dans le cas où l'invariant de Yamabe, défini par 
\begin{equation}
Y(M,C)=\inf_{g\in C}\frac{\int_M\scal_g\de v_g}{\Vol(M,g)^{\frac{n-2}2}},
\end{equation}
est strictement positif, ce qui est une restriction topologique et géométrique
assez importante, on a aussi l'identité classique pour $n\geq3$
\begin{equation}\label{dirac:eq2}
\inf_{g\in C}\lambda_1(M,\mathrm L_g)\Vol(M,g)^\frac2n=Y(M,C),
\end{equation}
qui nous dit que $\lambda^+(M,\mathrm L_g)=\lambda_1(M,\mathrm L_g)$ en 
donnant au passage une minoration optimale de $\lambda^+(M,\mathrm L_g)
\Vol(M,g)^\frac2n$, et qui fournit aussi en 
conjonction avec l'inégalité de Hijazi~(\ref{dirac:hijazi}) 
une minoration explicite de $\lambda_1(M,\chi,g)\Vol(M,g)^\frac2n$ pour 
$n\geq3$. C.~Bär a étendu dans \cite{ba92} cette inégalité 
à la dimension~2:
\begin{equation}
\lambda_1(S^2,\chi,g)\Vol(S^2,\chi,g)\geq4\pi.
\end{equation}
Pour les surfaces de genre supérieur ou égal à~1, l'invariant de Yamabe est 
négatif ou nul, et donc l'inégalité est triviale.

Signalons aussi que dans \cite{ah06b}, B.~Ammann et É.~Humbert étudient
l'existence de métriques minimisant $\lambda_2(M,\mathrm L_g)\Vol(M,g)^\frac2n$.

\subsection{Les invariants $\sigma$ et $\tau$}
On a vu que la majoration de la première valeur propre du laplacien dans 
une classe conforme permettait de définir un invariant différentiel
par un min-max. On peut définir des invariants similaires
à l'aide de l'opérateur de Dirac et du laplacien conforme et d'un max-min.

Un invariant de ce type a déjà été défini  
à l'aide de l'invariant de Yamabe. Cet invariant, introduit indépendamment
par O.~Kobayashi \cite{ko87} et R.~Schoen \cite{sc87}, est noté $\sigma(M)$ 
et appelé invariant de Schoen ou invariant de Yamabe différentiel:
\begin{equation}\label{dirac:sigma1}
\sigma(M)=\sup_gY(M,[g]),
\end{equation}
Comme pour toute classe conforme $C$ on a $Y(M,C)\leq n(n-1)\omega_n^{\frac2n}$
cette borne supérieure est finie, et $\sigma(M)\leq n(n-1)
\omega_n^{\frac2n}$.
Dans le cas ou $\sigma(M)$ est strictement positif, l'égalité~(\ref{dirac:eq2})
nous dit que c'est l'invariant annoncé pour le laplacien conforme:
\begin{equation}\label{dirac:sigma2}
\sigma(M)=\sup_{g}\inf_{\tilde g\in[g]}\lambda_1(M,\mathrm L_{\tilde g})
\Vol(M,\tilde g)^\frac2n.
\end{equation}
\begin{remarque}
L'invariant qu'on pourrait définir de la même manière avec 
$\lambda^+(M,\mathrm L_g)$
semble moins pertinent (et plus difficile à étudier) car le nombre
de valeurs propres négatives ou nulles de $\mathrm L$ dépend de la classe
conforme.
\end{remarque}

L'invariant $\sigma$ a été plus étudié que l'invariant $\nu$ de 
Friedlander-Nadirashvili vu au paragraphe~\ref{fonctions:scm} en raison de 
ses liens avec le problème de Yamabe, l'existence de métriques de courbure 
scalaire positive et la recherche de métriques d'Einstein 
(voir par exemple \cite{sc89}). Un certain
nombre de travaux récents portent sur l'évolution de $\sigma(M)$ quand
on procède à des chirurgies sur $M$ (voir notamment \cite{pe98}, 
\cite{pe00} et \cite{py99}). Je ne présenterai pas ici l'ensemble des
résultats obtenus, je me contenterai de rappeler quelques valeurs exactes
connues de l'invariant~$\sigma$.

En dimension quelconque, on sait calculer l'invariant~$\sigma$ sur les 
variétés suivantes (voir \cite{ko87} et \cite{sc89}) :
$$\sigma(S^n)=\sigma(\#i(S^{n-1}\times S^1))=n(n-1)\omega_n^{\frac2n},
\textrm{ et } \sigma(T^n\#N)=0$$
pour tout entier $i\in\N^*$ et toute variété $N$ telle que $\sigma(N)\geq0$.

Une attention particulière a été portée aux variétés de dimension~3:
pour tous entiers $i$, $j$, $k$ et $l$ tels que $i+j\geq1$, on a
$$\sigma(\#i(\R P^3)\#j(\R P^2\times S^1)\#k(S^2\times S^1)\#
l(S^2\ltimes S^1))=6\pi^{\frac23}=\frac{\sigma(S^3)}{2^{\frac23}},$$
où $S^2\ltimes S^1$ désigne le fibré non orientable en cercle sur $S^2$
(voir \cite{an07} et les références qui y sont données).
L'exemple de $\R P^3$, calculé par H.~L.~Bray et A.~Neves dans \cite{bn04},
vient conforter la conjecture selon laquelle pour les quotients de
la sphères, on a $\sigma(\Gamma\backslash S^n)=n(n-1)\left(\omega_n
/|\Gamma|\right)^{\frac2n}$.

 Dans \cite{ah06a}, B.~Ammann et É.~Humbert définissent l'invariant
correspondant pour l'opérateur de Dirac
\begin{equation}\label{dirac:tau}
\tau(M,\chi)=\sup_{g}\inf_{\tilde g\in[g]}\sqrt{\lambda_1(M,\chi,\tilde g)}
\Vol(M,\tilde g)^\frac1n,
\end{equation}
l'une des motivations étant que que l'inégalité de Hijazi le lie à 
l'invariant de Schoen:
\begin{equation}\label{dirac:tausig}
\tau(M,\chi)^2\geq\frac n{4(n-1)}\sigma(M).
\end{equation}
Dans ce cas aussi le nombre de valeurs propres nulles peut varier avec
la classe conforme, c'est la raison pour laquelle on considère 
$\lambda_1(M,\chi,\tilde g)$ et pas $\lambda^+(M,\chi,\tilde g)$.

On a la majoration $\tau(M,\chi)\leq\tau(S^n,\chi)=\frac n2\omega_n^{\frac1n}$
(voir \cite{am03}, \cite{ahm03} et \cite{ahm04}) et on connaît la valeur de 
$\tau$ sur quelques variétés : dans \cite{ah06a},
B.~Ammann et É.~Humbert déduisent de (\ref{dirac:tausig}) que
\begin{equation}
\tau(S^{n-1}\times S^1,\chi)=\frac n2\omega_n^{\frac1n}=\tau(S^n,\chi),
\end{equation}
et montrent que $\tau(T^2,\chi)=0$ ou $2\sqrt\pi$ selon la structure 
spinorielle $\chi$. Dans \cite{ah06c}, ils étendent ce dernier résultat aux 
autres surfaces orientables compactes et montrent en toute dimension que
$\tau(M,\chi)$ croit par adjonction d'anses. 

\subsection{Grandes valeurs propres}
Les valeurs propres non nulles de l'opérateur de Dirac et du laplacien
conforme ne sont pas bornées sur une classe conforme comme celle du laplacien
usuel :
\begin{theo}[\cite{aj07}]\label{dirac:sup}
Soit $M$ une variété compacte de dimension $n$ et $C$ une classe conforme
de métrique sur $M$. Si $n\geq3$, alors
$$\sup_{g\in C}\lambda^+(M,\mathrm L_g)\Vol(M,g)^\frac2n=+\infty,$$
et
$$\inf_{g\in C}\lambda^-(M,\mathrm L_g)\Vol(M,g)^\frac2n=-\infty.$$
Si $n\geq2$ et que $M$ est munie d'une structure spinorielle $\chi$, alors
$$\sup_{g\in C}\lambda^+(M,\chi,g)\Vol(M,g)^\frac2n=+\infty.$$
\end{theo}
La démonstration de ce théorème utilise des métriques, baptisées 
«~métriques Pinocchio~» dans \cite{ab00}, qui sont telles qu'un domaine
de la variété soit presque isométrique à la réunion d'un cylindre 
arbitrairement long et d'une demi-sphère collée à l'une de ses extrémité 
(le «~nez~» de la métrique). Pour l'opérateur de Dirac et en dimension 
$n\geq3$ on utilise la formule de Schrödinger-Lichnerovicz
\begin{equation}\label{dirac:sl}
\mathrm D^2=\nabla^*\nabla+\frac14\mathrm{scal}_g
\end{equation}
et le fait que la courbure scalaire est minorée sur le nez pour contrôler
la première valeur propre. Le laplacien conforme vérifie par définition
une relation semblable à (\ref{dirac:sl}), ce qui permet d'adapter facilement
la démonstration à cet opérateur. Le cas de la dimension~2 est plus technique 
car la courbure d'un cylindre est alors nulle, mais on peut s'aider du fait 
que toute métrique est localement conformément plate.

\section{Laplacien de Hodge-de~Rham}\label{formes}
\subsection{Définitions et notations}
 Le laplacien de Hodge-de~Rham agit sur l'espace $\Omega(M)$ des
formes différentielles de la variété $(M,g)$, et est défini par
$$\Delta=\de\codiff+\codiff\de,$$
où $\codiff$ désigne la codifférentielle. Son spectre en restriction
aux formes de degré $p$ sera noté
\begin{equation}
0=\lambda_{p,0}(M,g)<\lambda_{p,1}(M,g)\leq\lambda_{p,2}(M,g)\leq\ldots
\end{equation}
où les valeurs propres non nulles sont répétées s'il y a multiplicité.
La multiplicité de la valeur propre nulle, si elle existe, est un invariant
topologique: c'est le nombre de Betti $b_p(M)$. 

L'espace des $p$-formes coexactes est stable par le laplacien, et
on notera 
\begin{equation}
0<\mu_{p,1}(M,g)\leq\mu_{p,2}(M,g)\leq\ldots
\end{equation} 
le spectre du laplacien restreint à cet espace. La théorie de Hodge
nous dit que le spectre $(\lambda_{p,i}(M,g))_{i\geq1}$ est la réunion de 
$(\mu_{p,i}(M,g))_i$ et $(\mu_{p-1,i}(M,g))_i$. On a de plus
$\lambda_{p,i}(M,g)=\lambda_{n-p,i}(M,g)$ et donc 
$\mu_{p,i}(M,g)=\mu_{n-p-1,i}(M,g)$ car on considère des variétés sans bord.
Par conséquent, le spectre complet du laplacien se déduit des 
$\mu_{p,i}(M,g)$ pour $p\leq\frac{n-1}2$, et on se restreindra souvent
dans la suite à l'étude de ces valeurs propres. Rappelons aussi que 
$(\lambda_{p,0}(M,g))$ est le spectre du laplacien agissant sur les fonctions, 
et qu'il est aussi égal à ($\mu_{p,0}(M,g))$ si on excepte la valeur 
propre nulle.

 Le laplacien de Hodge-de~Rham n'est pas conformément covariant, mais
il possède une autre propriété très utile pour étudier les extrema
conformes du spectre, à savoir qu'on peut comparer les spectres de
deux métriques dont on connaît le rapport de quasi-isométrie :
\begin{proposition}[\cite{do82}]\label{formes:dodziuk}
Soit $g$ et $\tilde g$ deux métriques riemanniennes sur une variété 
compacte $M$ de dimension $n$, et $\tau$ une constante strictement positive. 
Si les deux métriques vérifient $\frac1\tau g\leq\tilde g\leq\tau g$, alors
$$\frac1{\tau^{3n-1}}\lambda_{p,k}(M,g)\leq\lambda_{p,k}(M,\tilde g)
\leq\tau^{3n-1}\lambda_{p,k}(M,g),$$
pour tout entiers $k\geq0$ et $p\in[0,n]$.
\end{proposition}
 Si deux métriques vérifient $\frac1\tau g\leq\tilde g\leq\tau g$, on a
aussi $\frac1\tau h^2g\leq h^2\tilde g\leq\tau h^2g$ pour toute fonction
$h$ strictement positive. Par conséquent, pour montrer qu'on peut faire 
tendre une valeur propre vers $0$ ou $+\infty$ dans une classe conforme,
il suffit de le montrer pour la classe conforme d'une métrique $g$ bien
choisie, le résultat se transposant automatiquement à n'importe quelle
autre classe conforme. 
\begin{remarque}
Outre la continuité du spectre pour la topologie $C^0$, la 
proposition~\ref{formes:dodziuk} 
implique que la dimension du noyau du laplacien ne dépend pas de la métrique.
On ne peut donc pas espérer un tel résultat pour l'opérateur de Dirac.
\end{remarque}
\subsection{Prescription du spectre}
Le comportement du spectre du laplacien de Hodge-de~Rham dans une classe
conforme est en partie analogue à celui de l'opérateur de Dirac, mais on va 
voir que pour certains degrés toute rigidité disparaît et qu'on peut
prescrire arbitrairement le début spectre.

Le premier résultat qui a été obtenu est que, comme pour l'opérateur
de Dirac, on peut rendre les valeurs propres arbitrairement grande~:
\begin{theo}[\cite{ces06}]\label{formes:gvp}
Si $M$ est une variété riemannienne compacte de dimension $n\geq3$ et $C$ 
une classe conforme sur $M$, alors
\begin{equation}
\sup_{g\in C}\inf_{1\leq p\leq\frac n2}\mu_{p,1}(M,g)\Vol(M,g)^{\frac2n}=
+\infty.
\end{equation}
\end{theo}
\begin{remarque} Ce théorème est énoncé dans \cite{ces06} pour les
variétés de dimension $n\geq4$, mais la démonstration s'applique
parfaitement aux 1-formes coexactes en dimension~3.
\end{remarque}
Les métriques utilisée par B.~Colbois et A.~El~Soufi pour
démontrer le théorème~\ref{formes:gvp} sont des «~métriques
Pinocchio~», comme dans le cas de l'opérateur de Dirac, mais
les outils d'analyse sont différents~: on utilise le lemme de 
McGowan (voir \cite{mc93}, lemme~2.3), qui permet de minorer le spectre 
de la variété en fonction des spectres d'une famille de domaines qui 
recouvre $M$. Le théorème~\ref{formes:gvp} n'est bien sûr pas valable
pour $p=0$ comme on l'a déjà vu.

Ce résultat appelle  naturellement la question de savoir si l'existence
de petites valeurs propres dans une classe conforme pour spectre des 
fonctions se généralise aux formes~:
\begin{question}[\cite{co04}]\label{formes:q1}
Étant donnés $1\leq p\leq\frac n2$ et $k\geq1$, a-t-on
\begin{equation}
\inf_{g\in C}\mu_{p,k}(M,g)\Vol(M,g)^{\frac2n}=0\ ?
\end{equation}
\end{question}

Dans le cas du laplacien de Hodge-de~Rham, cette question a une résonance 
particulière. En effet, le noyau du laplacien étant canoniquement isomorphe
à la cohomologie de de~Rham, chercher les petites valeurs propres du 
laplacien revient à déterminer quelles déformations de la métrique
tendent à «~créer~» de la cohomologie. Ce problème a pour l'instant
été étudié dans deux situations géométriques particulières~: les limites 
adiabatiques associées aux feuilletages riemanniens (voir par exemple 
\cite{ALK00} et les références qui y sont données) et les effondrement à 
courbure bornée (voir \cite{Ja05} pour une présentation synthétique
des résultats à ce sujet). Ces deux situations se recoupent partiellement
---~sans que l'on puisse réduire l'une à l'autre~---, mais sont
très différentes des déformations à volume et classe conforme fixée.
Une déformation conforme est ponctuellement isotrope mais globalement
inhomogène, alors que les limites adiabatiques sont des déformations 
homogènes et anisotropes (en tout point, on décompose l'espace tangent en 
la somme $T_xM=H\oplus V$ de deux espaces muni chacun d'une métrique $g_H$ 
et $g_V$, et on considère la famille de métrique 
$g_\varepsilon=g_H\oplus\varepsilon^2g_V$ avec $\varepsilon\to0$), et 
que les effondrements à courbure bornée permettent, sous certaines 
conditions topologiques, une anisotropie encore plus grande.

J'ai apporté une première partie de la réponse à la question \ref{formes:q1}
dans \cite{ja06c} en montrant que pour presque tous les degrés, on
peut faire tendre un nombre arbitraire de valeur propres vers zéro :
\begin{theo}\label{formes:pvp}
Soit $M$ une variété compacte sans bord de dimension $n\geq5$, $k$ l'entier 
tel que $n=2k+3$ ou $n=2k+4$ et $C$ une classe conforme
sur $M$. Pour tout réel $V>0$ et toute suite d'entiers positifs ou nuls 
$N_1,N_2,\ldots,N_{k}$, 
il existe une famille de métriques $(g_\varepsilon)_{0<\varepsilon<1}$ 
contenue dans $C$ et une constante $c>0$ telles que 
$\mu_{p,N_p}(M,g_\varepsilon)<\varepsilon$ et 
$\mu_{p,N_p+1}(M,g_\varepsilon)>c$ pour tout $1\leq p\leq k$, 
$\mu_{k+1,1}(M,g_\varepsilon)>c$ et $\Vol(M,g_\varepsilon)=V$.
\end{theo}
On peut donc prescrire le nombre de petite valeurs propres pour les
formes coexactes. Les degrés $\frac n2$ et $\frac n2-1$ si $n$ est pair,
et le degré $[\frac n2]$ si $n$ est impair font exception, on verra
au paragraphe \ref{formes:min} qu'on ne peut pas faire tendre les
valeurs propres vers zéro dans ces cas.

 La construction géométrique intervenant dans la démonstration du théorème 
\ref{formes:pvp} consiste à écraser la métrique en dehors du voisinage
tubulaire d'une ou plusieurs sous-variétés, les formes
harmoniques de ces sous-variétés peuvent alors s'étendre en des formes
test ayant un petit quotient de Rayleigh. La difficulté est alors de 
contrôler précisément le nombre de petites valeurs propres. Cette technique 
échoue pour les degrés proches de $\frac n2$ car la norme $L^2$ des 
formes devient alors conformément invariante ou presque invariante.

L'existence de de grandes et de petites valeurs propres à volume
fixé pour le laplacien de Hodge-de~Rham a permis à P.~Guérini de
montrer dans \cite{gu04} qu'on pouvait prescrire le volume et
le début du spectre sur une variété compacte quelconque, généralisant
ainsi aux formes différentielles les résultats obtenus par Y.~Colin de
Verdière et J.~Lohkamp pour les fonctions. Les théorèmes \ref{formes:gvp} 
et \ref{formes:pvp} permettent de prescrire simultanément le volume, le 
début du spectre et la classe conforme~:
\begin{theo}[\cite{ja06c}]\label{formes:presc}
Soit $M$ une variété compacte, connexe, orientable et sans bord de dimension 
$n=2k+3$ ou $2k+4$ où $k\in\N^*$,
$C$ une classe conforme de métriques riemanniennes sur $M$,
$V$ un réel strictement positif et $N\geq1$ un entier. 
On se donne pour 
tout entier $p\in\{1,\ldots,k\}$ une suite de réels 
$0<\nu_{p,1}<\nu_{p,2}<\ldots<\nu_{p,N}$.

Il existe une métrique $g\in C$ telle que
\begin{itemize}
\item $\mu_{p,i}(M,g)=\nu_{p,i}$ pour tout $i\leq N$ et $p\in\{1,\ldots,k\}$;
\item $\mu_{k+1,1}(M,g)>\sup_{p,i}\{\nu_{p,i}\}$;
\item $\Vol(M,g)=V$.
\end{itemize}
\end{theo}

\begin{remarque}
La minoration $\mu_{k+1,1}(M,g)>\sup_{p,i}\{\nu_{p,i}\}$ assure qu'on a
l'égalité
$\lambda_{k+1,i}(M,g)=\mu_{k,i}(M,g)$ pour $i\leq N$. On peut donc
prescrire les $N$ premières valeurs propres des $(k+1)$-formes, les formes 
propres correspondantes étant alors exactes, de valeurs propres égales
à $(\mu_{k,i}(M,g))_{i=1}^N$. Si $n$ est impair,
on prescrit ainsi le spectre en tout degré $2\leq p\leq n-2$. En
dimension paire, le degré $p=n/2=k+2$ fait exception. En degré $1$ et $n-1$
on ne prescrit pas arbitrairement le début du spectre car on ne contrôle 
pas les $\mu_{0,i}(M,g)$, mais on peut assurer que les valeurs 
$\nu_{1,1},\ldots,\nu_{1,N}$ sont contenues dans 
$(\lambda_{1,i}(M,g))_{i\geq1}$ et $(\lambda_{n-1,i}(M,g))_{i\geq1}$.
\end{remarque}
Le théorème \ref{formes:presc} contraste avec la rigidité apparaissant 
dans les cas du laplacien agissant sur les fonctions et de l'opérateur 
de Dirac. Dans le cas du laplacien de Hodge-de~Rham il reste une certaine 
souplesse, mais de justesse puisque certains degrés font exception.

 On peut se demander si la prescription du spectre en degré $n/2$ 
serait possible si on supprimait la contrainte sur le volume de la variété,
ou, de manière équivalente, si on peut faire tendre une unique valeur
propre vers zéro sans fixer le volume. Ce problème reste à éclaircir (voir
question~\ref{formes:q}).

\subsection{Inégalités de Sobolev pour les formes différentielles}

Pour achever de répondre à la question \ref{formes:q1} nous aurons à
faire appel à des inégalités de Sobolev pour les formes différentielles,
dont certaines constantes s'avèrent être des invariants conformes, et 
que nous allons rappeler ici. Elles reposent essentiellement
sur le fait que la différentielle extérieure définit un complexe
elliptique sur l'espaces des formes différentielles. Les ouvrages
consacrés aux opérateurs (pseudo)différentiels et aux inégalités
elliptiques développent rarement le cas des opérateurs agissant
sur les sections d'un fibré vectoriel ; on peut toutefois se référer
au chapitre~6 de \cite{mo66}, et une définition précise d'un 
complexe elliptique est donnée dans \cite{ra05} (chapitre~9) et 
dans \cite{ag94}.

 On se donne une variété riemannienne compacte $(M^n,g)$, et deux réels
$r,s>1$ tels que $\frac1s-\frac1r\leq\frac1n$. L'inégalité de Sobolev
classique sur les fonctions assure l'existence de deux constantes $A,B>0$
telles que pour toute fonction $f\in C^\infty(M)$,
\begin{equation}
\|f\|_{L^r}\leq A\|\de f\|_{L^s}+B\|f\|_{L^1}.
\end{equation}
Cette inégalité se généralise aux formes différentielles en faisant
intervenir la codifférentielle :
\begin{theo}\label{formes:ineg1}
Il existe trois constantes $A,B,C>0$ dépendant de $g$, $r$ et $s$ telles 
que pour toute forme différentielle $\omega\in\Omega(M)$, on a
\begin{equation}\label{formes:ineg2}
\|\omega\|_r\leq A\|\de\omega\|_s+B\|\codiff\omega\|_s+C\|\omega\|_1.
\end{equation}
\end{theo}
Le théorème \ref{formes:ineg1} est un corollaire immédiat de l'inégalité
elliptique associée à l'opérateur $(\de+\codiff)$ 
et de la continuité de l'injection de Sobolev.

Nous aurons en fait besoin d'une inégalité prenant une forme différente:
\begin{theo}\label{formes:theq}
Soit $(M,g)$ une variété compacte de dimension $n$, et deux réels $r,s>1$
tels que $\frac1s-\frac1r\leq\frac1n$.
Il existe une constante $c(M,g,r,s)>0$ telle que pour toute $p$-forme 
$\omega\in\Omega^p(M,g)$, on a
\begin{equation}\label{formes:eq2}
\inf_{\de\varphi=0}\|\omega-\varphi\|_r\leq c\|\de\omega\|_s.
\end{equation} 
\end{theo}
On peut démontrer ce théorème en utilisant le fait qu'en restriction à 
l'orthogonal des formes harmoniques le dernier terme de 
l'inégalité~(\ref{formes:ineg2}) est superflu (quitte à modifier les 
constantes $A$ et $B$), et en appliquant cette 
inégalité à la projection orthogonale de $\omega$ sur les formes coexactes.
On trouve aussi dans l'appendice de \cite{gt06} une démonstration 
détaillée du théorème~\ref{formes:theq} se basant sur l'inégalité
elliptique du laplacien de Hodge-de~Rham.

 L'inégalité du théorème~\ref{formes:theq} nous intéressera pour une 
particularité apparaissant quand on fixe le degré $p$ et qu'on pose 
$r=\frac np$ et $s=\frac n{p+1}$. En effet, les normes $\|\omega-\varphi\|_r$ 
et $\|\de\omega\|_s$ sont alors conformément invariantes, et par conséquent 
la constante optimale dans l'inégalité (\ref{formes:eq2}) en restriction 
aux $p$-formes est un invariant conforme, qu'on notera
\begin{equation}\label{formes:eq3}
K_p(M,[g])=\sup_{\omega\in\Omega^p(M)}\inf_{\de\varphi=0}
\frac{\|\omega-\varphi\|_{\frac np}}{\|\de\omega\|_{\frac n{p+1}}},
\end{equation}

\begin{remarque}
On ne connaît actuellement aucune estimée de la constante $K_p$.
\end{remarque}

\subsection{Minoration du spectre}\label{formes:min}

Sur les variétés de dimension $n\geq3$, la constante de Sobolev $K_p$ que 
nous avons défini au paragraphe précédent permet de minorer le spectre 
du laplacien de Hodge-de~Rham dans les cas qui ne sont pas couverts par le 
théorème \ref{formes:pvp}~:
\begin{theo}[\cite{ja07}]\label{formes:thmin}
Soit $M^n$ une variété compacte de dimension $n\geq3$, $C$ une classe
conforme de métriques sur $M$. Pour
toute métrique $g\in C$, on a
$$\mu_{\left[\frac n2\right],1}(M,g)\Vol(M,g)^{\frac2n}\geq 
K_{\left[\frac n2\right]}(M,C)^{-2}.$$
Si $n$ est pair, cette inégalité est optimale. 
\end{theo}
\begin{remarque} Quand $n$ est pair, on a $\mu_{\frac n2-1,i}(M,g)=
\mu_{\frac n2,i}(M,g)$. Les théorèmes \ref{formes:pvp} et \ref{formes:thmin}
répondent donc complètement à la question \ref{formes:q1}.
\end{remarque}
\begin{remarque}
Une conséquence du théorème~\ref{formes:thmin} est que, paradoxalement, 
les déformations de la métriques qui apparaissent dans la démonstration
du théorème~\ref{formes:pvp} et qui tendent à créer de l'homologie 
ne produisent pas de petites valeurs propres en degré médian.
\end{remarque}
\begin{remarque}
On peut facilement trouver des variétés qui s'effondrent en faisant tendre 
$\mu_{\left[\frac n2\right],1}(M,g)\Vol(M,g)^{\frac2n}$ vers zéro. On
ne peut donc pas donner de majoration uniforme de 
$K_{\left[\frac n2\right]}(M,C)$ en fonction de bornes sur le diamètre et 
la courbure de la variété pour une métrique $g\in C$.
\end{remarque}
\begin{remarque}
Une conséquence assez inattendue du théorème~\ref{formes:thmin} est
la construction, pour tout $n\geq4$ et tout $1\leq p\leq\frac n2$ 
de variétés $M$ de dimension $n$  telles que la multiplicité de
$\mu_{p,1}(M,g)$ soit arbitrairement grande (voir \cite{ja06a}).
\end{remarque}

 On peut se demander s'il existe d'autres rigidités, par exemple:
\begin{question}\label{formes:q}
Le rapport $\mu_{\left[\frac n2\right],k+1}(M,g)/
\mu_{\left[\frac n2\right],k}(M,g)$ est-t-il borné sur une classe conforme ?
\end{question}

Comme pour les fonctions, on peut définir un spectre conforme pour
les formes différentielles par
\begin{equation}
\mu_k^c(M,C)=\inf_{g\in C}\mu_{\left[\frac n2\right],k}(M,g)
\Vol(M,g)^{\frac2n}.
\end{equation}
Le théorème \ref{formes:thmin} donne la valeur de la première
valeur propre conforme en fonction de la constante $K_p$ en dimension paire,
et une minoration en dimension impaire. On peut aussi donner une majoration
du spectre conforme en fonction de celui de la sphère~:
\begin{theo}[\cite{ja07}]\label{formes:minmaj}
Pour toute classe conforme $C$ sur la variété compacte $M$ et tout $k\geq1$, 
on a
$$\mu_k^c(M,C)\leq\mu_k^c(S^n,C_\textrm{\emph{can}})\leq k^{\frac2n}
\mu_1^c(S^n,C_\textrm{\emph{can}}),$$
où $C_\textrm{\emph{can}}$ est la classe conforme de la métrique canonique 
de la sphère.
\end{theo}
 
Un problème naturel est de déterminer si la borne inférieure $\mu_k^c(M,C)$
est atteinte par une métrique régulière, en particulier si $k=1$. Dans
\cite{ja07}, on établit le critère suivant~:
\begin{theo}\label{formes:mincn}
Soit $g$ une métrique lisse de volume~1 sur $M$ telle que l'on ait 
$\mu_{\left[\frac n2\right],1}(M,g)=\mu_1^c(M,[g])$. Si $n=3\textrm{ mod }4$,
alors il existe une $\left[\frac n2\right]$-forme propre coexacte non nulle
de valeur propre $\mu_1^c(M,[g])$ et
de longueur constante. Si $n$ est pair, toutes les formes propres coexacte 
de degré $\frac n2-1$ et de valeur propre $\mu_1^c(M,[g])$ sont de 
longueur constante.
\end{theo}
\begin{remarque}\label{formes:minrq} Sur les sphères munies de leur 
métrique canonique, 
les premières formes propres ne sont pas de longueur constantes (voir
l'appendice de \cite{gm75}). Contre toute attente, la métrique canonique
de la sphère n'est donc pas extrémale pour la première valeur propre !
Le théorème~\ref{formes:mincn}
laisse présager que les métriques extrémales pour les formes différentielles
sont en général assez différentes de celles des fonctions.
\end{remarque}
En dimension~4, on peut déduire du théorème \ref{formes:mincn} une
condition nécessaire sur la topologie pour l'existence de métriques
extrémales lisses~:
\begin{cor}\label{formes:mincor}
 Si $M$ est une variété compacte de dimension~4 et de
caractéristique d'Euler non nulle, il n'existe pas de métrique régulière
$g$ de volume~1 telle que $\mu_{1,1}(M,g)=\mu_1^c(M,[g])$.
\end{cor}
Dans le cas des fonctions on connaît des conditions nécessaires sur
les métriques extrémales (cf.~\cite{esi03}), mais pas d'obstruction 
topologique comme celle du corollaire \ref{formes:mincor}.

On ne connaît finalement aucun exemple de métrique extrémale pour ce 
problème, ni même de valeur explicite de $\mu_1^c(M,[g])$. Il existe
toutefois des candidats sérieux, comme certaines variétés effondrées
ayant une petite valeur propre dont la forme propre est de longueur
constante, en particulier dans les situations topologiques les plus simples,
par exemple:
\begin{itemize}
\item les fibrés en cercles sur les
surfaces munies d'une métrique effondrée adaptée (au sens de \cite{cc00},
définition 2.1) ;
\item les sphères munies d'une métrique de Berger ; 
\item les fibrés en tores sur le cercle munis d'une métrique homogène 
effondrée produisant une petite valeur propre (voir \cite{ja03}).
\end{itemize}

On peut utiliser le spectre conforme des formes différentielles pour 
définir des invariants différentiels de la variété comme on l'a vu 
pour les autres opérateurs, en posant
\begin{equation}
\mu_k(M)=\sup_g\mu_k^c(M,[g]).
\end{equation}
Le théorème \ref{formes:minmaj} assure que cette borne supérieure est
finie. On ne connaît cependant la valeur de ces invariants sur aucune
variété, deux raisons en étant que $\mu_k(M)$ n'est pas défini en
dimension~2, contrairement à ce qui se passe pour le laplacien sur
les fonctions ou l'opérateur de Dirac, et que le laplacien de Hodge-de~Rham
n'est pas conformément covariant. Même sur les sphères, on
peut seulement affirmer que $\mu_k(S^n)=\mu_k^c(S^n,C_\textrm{can})$
sans en expliciter la valeur (voir le théorème \ref{formes:mincn} et
la remarque \ref{formes:minrq} ci-dessus). La question de la trivialité
des ces invariants reste ouverte.

\noindent Pierre \textsc{Jammes}\\
Université d'Avignon et des pays de Vaucluse\\
Laboratoire d'analyse non linéaire et géométrie (EA 2151)\\
33 rue Louis Pasteur\\
F-84000 Avignon\\
\texttt{Pierre.Jammes@univ-avignon.fr}
\end{document}